\documentclass{article}

\makeatletter{}\usepackage[boxed]			{algorithm2e}
\usepackage			{amsmath}
\usepackage			{amsfonts}
\usepackage			{array}
\usepackage			{comment}
\usepackage			{enumerate}
\usepackage[mathscr]		{eucal}
\usepackage			{graphicx}
\usepackage[utf8]		{inputenc}
\usepackage			{latexsym}
\usepackage			{longtable}
\usepackage			{multirow}
\usepackage			{rcs}
\usepackage			{theorem}
\usepackage			{xspace}

\newtheorem{theorem}{Theorem}

\newcommand{\EndProof}{\rule{2ex}{2ex}}
\newcommand{\uscy}{upper semi-continuity\xspace}
\newcommand{\usc}{upper semi-continuous\xspace}

\newcommand{\bleq}[2]{\begin{equation}\label{#1}{#2}\end{equation}}

\DeclareMathOperator	{\cl}	{{\mathrm cl}}

\DeclareMathOperator	{\co}	{{\mathrm co}}

\DeclareMathOperator	{\intr}	{{\mathrm int}}

\DeclareMathOperator	{\dom}	{{\mathrm dom}}

\newcommand{\R}{\mathbb{R}}
\newcommand{\xb}{\bar x}

\DeclareMathOperator*{\epi}{epi}

\newcommand{\xs}{{x^\star}}

\newcommand{\enProof}{{\bf Proof}\xspace}
\newcommand{\cC}{{\cal{C}}}

\newtheorem{enlem}{Lemma}
\newtheorem{definition}{Definition}

\title{Finite-Value Superiorization for Variational Inequality Problems}
\author{E.A. Nurminski}
\date{October 12, 2016}
\begin{document}
\maketitle
\begin{abstract}
\makeatletter{}The main goal of this paper is to present the application
of a superiorization methodology
to solution of variational inequalities.
Within this framework
a variational inequality operator is considered as a small perturbation
of a convex feasibility solver
what allows to construct a simple iteration algorithm.
The specific features of variational inequality problems
allow to use a finite-value perturbation which may be
advantageous from computational point of view.
The price for simplicity and finite-value is that
the algorithm provides an approximate solution of
variational inequality problem with a prescribed
coordinate accuracy.
 
\end{abstract}
\section*{Introduction}
\makeatletter{}\label{intro}
This paper presents a variant of superiorization methodology
for solution of variational inequalities (VI) problems.
VI became one of the common tools for representing
many problems in physics, engineering, economics, computational biology,
computerized medicine, to name but a few,
which extend beyond optimization, see
\cite{intro03}
for the extensive review of the subject.
Apart from the mathematical problems connected with the characterization
of solutions and development of the appropriate algorithmic tools to find them,
modern problems offer significant implementation challenges
due to their nonlinearity and large scale.
It leaves just a few options for the algorithms development
as it occurs in the others related fields like convex feasibility (CF) problems
\cite{bb96} as well.
One of these options is to use
fixed point iteration methods
with various attraction properties toward the solutions,
which have low memory requirements and simple
and easily parallelized iterations.
These schemes are quite popular for convex optimization and CF problems
but they need certain modifications to be applied to VI problems.
The idea of modification can be related to some approaches put forward for
convex optimization and CF problems
\cite{nurmiDAN08,nurmiCMMP08,sup2010}
and which is becoming known as superiorization technique (see also \cite{psm-su} for the general description).

From the point of view of this approach the optimization problem
\bleq{op}{\min f(x),\quad x \in X}
or VI problem to find \(\xs \in X\) such that
\bleq{vi}{ F(\xs)(x - \xs) \geq 0, \quad x \in X }
are conceptually divided into the feasibility problem \(x \in X \) and
the second-stage optimization or VI problems next.
Then we may consider these tasks to a certain extent separately and make use of their specifics
to apply the most suitable algorithms for feasibility and optimization/VI parts.
The problem is to combine these algorithms in a way which provides
the solution of the original problems (\ref{op}) or (\ref{vi}).
As it turns out these two tasks can be merged together
under rather reasonable conditions
which basically require
a feasibility algorithm to be resilient with respect to diminishing perturbations
and the second-stage algorithm to be globally convergent
over the feasible set or its expansions.
In the field of optimization this idea, as the author may say, was applied
quite early
in \cite{novosib74}
to non-stationary
extremum problems of the kind
\( \min_x f(x, \tau) = f(\xs(\tau), \tau) \) where parameter(s) \(\tau\) is changing in a certain
unknown in advance way.
The problem of interest for instance for the design of real-time systems
is to track the solution \(\xs(\tau)\) with hopefully asymptotic convergence
to the trajectory of these minima.
Again it is quite clear that the algorithm for such problems must be resilient
with respect to changes in \(f(x,\tau)\) which can be considered as perturbations.
The early results \cite{novosib74} demonstrated basically that the gradient technique satisfies these
requirements provided that step-sizes of the gradient process diminish slower
than changes in \(\tau\).
Further developments of superiorization methodology may widen the choice
of algorithms for this task.

The theory of the above mentioned approaches was developed mainly for optimization problems,
where asymptotic convergence
of these methods was studied when the number of iterations tends to infinity.
Needless to say that in practice it is always finite and depends on stopping criteria
of different nature and the starting position,
take text-book Newton method as an example.
Then one may see that the stopping point of the previous run, even if
it does not provide a solution for the problem at hand,
may significantly improve computational characteristics of the following runs
even if they are using the same algorithm but with the different balance between
feasibility and optimality goals.
Actually we may not even insist on exact solution in favor of decomposition and parallelization,
provided that the deviation from the solution is under control.
It is exactly what happened here where the superiorized feasibility algorithm is able to guarantee
only
a finite accuracy of a solution of VI problems.

The established asymptotic convergence for each algorithm used at each run is nevertheless essential
to guarantee the attainability of a stopping criteria or a condition of sufficient improvements.
Study of the rate of asymptotic convergence of different methods may reveal essential characteristics
of the problem which influence the computational efficiency and which therefore should be taken care of
as intermediate goals for different stages of solution.
This article is devoted to the idea of using a finite-values perturbation which hopefully
takes away one of the requirements which slows down the convergence of iteration algorithms.
In the results obtained so far
for
the feasibility and optimality algorithms 
to be successfully combined together we had to ensure that the second-stage optimizing algorithm
produces small and diminishing steps which are the perturbations for the feasibility algorithm.
In the native superiorization technology even the summability of the perturbations is assumed
and considered as a distinctive feature \cite{bdhk2007,psm-su}.
In the similar developments \cite{nurmiDAN08,nurmiCMMP08} this requirement is lifted,
however
it is still necessary to have vanishing optimization steps.
This slows down the convergence to overall solution of (\ref{op})
in the same way as it does in the penalty function method
which consists in the solution of the auxiliary problem
of the kind
\bleq{penty}{ \min_x \{\,\Phi_X(x) + \epsilon f(x)\,\}  = \Phi_X(x_\epsilon) + \epsilon f(x_\epsilon)}
where \(\Phi_X(x) = 0 \) for \(x \in X\) and \(\Phi_X(x) > 0 \) otherwise.
The term \(\epsilon f(x)\) can be considered as the perturbation of the feasibility
problem \( \min_x \Phi(x)\)
and for classical smooth penalty functions
the penalty parameter \(\epsilon > 0\) must tend to zero to guarantee convergence of
\(x_\epsilon\) to the solution of (\ref{op}) as it takes place in the superiorization theory.
Therefore such situation may be called infinitesimal superiorization.
Definitely it makes the objective function \(f(x)\)
less influential in solution process of (\ref{op}) and hinders the optimization.

To overcome this problem the exact penalty functions
\(\Psi_X(\cdot)\)
can be used which provide the exact solution of (\ref{op})
\[
\min_x \{\,\Psi_X(x) + \epsilon f(x)\,\}  = \epsilon f(\xs)
\]
for \(\epsilon > 0\) and small enough under rather mild conditions.
In the spirit of what was said above it may be called
the finite-value superiorization.
The price for the conceptual simplification of the solution of
(\ref{vi})
is the inevitable non-differentiability
of the penalty function \(\Psi_X(x)\) and the corresponding worsening of
convergence rates for instance for gradient-like methods
(see \cite{nest83,nest04} for comparison).
Nevertheless the idea has a certain appeal, keeping in mind successes
of nondiffereniable optimization, and
the similar approaches with necessary modifications were used for VI problems starting from
\cite{browder66}
and followed by
\cite{konn2003,ansil2010,koku2011,micpat2013}
among others.

Here we introduce a geometrical notion of a sharp penalty mapping
for which it is possible to prove the existence
of a finite penalty constant which allows to suggest for monotone VI problem
the iteration algorithm with the operator which is strongly oriented with respect
to the solution of VI outside the given neighborhood of solution.
It allows also to get rid of the Slater condition for the functional constraint which simplifies
a theory for solving VI problems with equality constraints.

Next we prove an approximate convergence of the iteration method
with the penalized variational operator
where some care should be taken to keep the iteration process bounded.
Toward this we use something like restarts from a point
inside the area of possible locations of a solution bounded by a certain large enough ball.
It is proved that in this case only finite number of restarts requires
and convergence of the algorithm follows
from certain general conditions for iteration processes
developed early in \cite{nurmi-rus91}. 
\section{Notations and preliminaries}
\makeatletter{}\label{nots}
Let \(E\) denotes a finite-dimensional space with the inner product \(xy\) for \(x, y \in E\),
and the standard Euclidean norm \(\|x\| = \sqrt{xx}\).
The one-dimensional \(E\) is denoted nevertheless as \(\R\) and \(\R_\infty = \R \cup \{\infty\}.\)
The unit ball in \(E\) will be denoted as \(B = \{x: \|x\| \leq 1 \}\).
The space of bounded closed convex subsets of \(E\) is denoted as \(\cC(E)\).
For any \(X \subset E\) we denote its interior as \(\intr(X)\).
The closure of a set \(X\) is denoted as \(\cl\{X\}\) and its boundary as \(\partial X\).
The distance function \(\rho(x,X)\) between point \(x\) and set \(X \subset E\)
is defined as
\(\rho(x,X) = \inf_{y \in X} \| x - y \|\).
The norm of a set \(X\) is defined as \( \|X\| = \sup_{x \in X} \|x\| \).

The sum of two subsets \(A\) and \(B\) of \(E\) is denoted as \(A + B\) and
understood as \(A+B = \{ a + b, a \in A, b \in B \}\).
If \(A\) is a singleton \(\{a\}\) we write just \( a + B \).

Any open subset of \(E\) containing zero vector is called a neighborhood of zero in \(E\).
We use the standard definition of \uscy of set-valued mappings:
\begin{definition}
A  set-valued mapping \(F: E \to \cC(E)\)
is called \usc
if at any point \(\bar x\) for any neighborhood of zero \(U\)
there exists a neighborhood of zero \(V\) such that \( F(x) \subset F(\bar x) + U \)
for all \(x \in \bar x + V\).
\end{definition}
From the point of view of VI the most studied class of set-valued mappings is probably the monotone ones.
\begin{definition}
A  set-valued mapping \(F: E \to \cC(E)\) is called a monotone
if \( (f_x - f_y)(x - y ) \geq 0 \) for any \(x, y  \in E\)
and \(f_x \in F(x), f_y \in F(y)\).
\end{definition}
We use standard notations of convex analysis: if \(h:E \to \R_\infty\) is a convex function,
then \(\dom(h) = \{ x: h(x) < \infty \}\),
\( \epi h = \{ (\mu, x) : \mu \geq h(x), x \in \dom(h) \} \subset \R \times E \),
the sub-differential of \(h\) is defined as follows:
\begin{definition}
For a convex function \(h: E \to \R_\infty\)
a sub-differential of \(h\) at point \(\bar x \in \dom(h)\) is 
the set \(\partial h(\bar x)\) of vectors \(g\) such that
\( h(x) - h(\bar x) \geq g(x - \bar x) \) for any \(x \in \dom(h)\).
\end{definition}
This defines
a convex-valued upper semi-continuous maximal monotone set-valued mapping
\(\partial h: \intr(\dom(h)) \to \cC(E)\).
At the boundaries of \(\dom(h)\) the sub-differential of \(h\) may or may not exists.
For differentiable \(h(x)\) the classical gradient of \(h\) is denoted as \(h'(x)\).

We define the convex envelope of \(X \subset E\) as follows.
\begin{definition}
An inclusion-minimal set \(Y \in \cC(E)\) such that \(X \subset Y\)
is called a convex envelope of \(X\)
and denoted as \(\co(X)\).
\end{definition}

Our main interest is in finding a solution \(\xs\) of a following finite-dimensional VI problem
with a single-valued operator \(F(x)\):
\bleq{vip}{ \mbox{Find } \xs \in X \subset \cC(E)
\mbox{ such that }
F(\xs)(x - \xs) \geq 0 \mbox{ for all } x \in X.}
This problem has its roots in convex optimization and for \(F(x) = f'(x) \) VI (\ref{vip}) is
the geometrical formalization of the optimality conditions for (\ref{op}).

\begin{figure}
\begin{center}
\makeatletter{}\begingroup
  \makeatletter
  \providecommand\color[2][]{    \GenericError{(gnuplot) \space\space\space\@spaces}{      Package color not loaded in conjunction with
      terminal option `colourtext'    }{See the gnuplot documentation for explanation.    }{Either use 'blacktext' in gnuplot or load the package
      color.sty in LaTeX.}    \renewcommand\color[2][]{}  }  \providecommand\includegraphics[2][]{    \GenericError{(gnuplot) \space\space\space\@spaces}{      Package graphicx or graphics not loaded    }{See the gnuplot documentation for explanation.    }{The gnuplot epslatex terminal needs graphicx.sty or graphics.sty.}    \renewcommand\includegraphics[2][]{}  }  \providecommand\rotatebox[2]{#2}  \@ifundefined{ifGPcolor}{    \newif\ifGPcolor
    \GPcolorfalse
  }{}  \@ifundefined{ifGPblacktext}{    \newif\ifGPblacktext
    \GPblacktexttrue
  }{}    \let\gplgaddtomacro\g@addto@macro
    \gdef\gplbacktext{}  \gdef\gplfronttext{}  \makeatother
  \ifGPblacktext
        \def\colorrgb#1{}    \def\colorgray#1{}  \else
        \ifGPcolor
      \def\colorrgb#1{\color[rgb]{#1}}      \def\colorgray#1{\color[gray]{#1}}      \expandafter\def\csname LTw\endcsname{\color{white}}      \expandafter\def\csname LTb\endcsname{\color{black}}      \expandafter\def\csname LTa\endcsname{\color{black}}      \expandafter\def\csname LT0\endcsname{\color[rgb]{1,0,0}}      \expandafter\def\csname LT1\endcsname{\color[rgb]{0,1,0}}      \expandafter\def\csname LT2\endcsname{\color[rgb]{0,0,1}}      \expandafter\def\csname LT3\endcsname{\color[rgb]{1,0,1}}      \expandafter\def\csname LT4\endcsname{\color[rgb]{0,1,1}}      \expandafter\def\csname LT5\endcsname{\color[rgb]{1,1,0}}      \expandafter\def\csname LT6\endcsname{\color[rgb]{0,0,0}}      \expandafter\def\csname LT7\endcsname{\color[rgb]{1,0.3,0}}      \expandafter\def\csname LT8\endcsname{\color[rgb]{0.5,0.5,0.5}}    \else
            \def\colorrgb#1{\color{black}}      \def\colorgray#1{\color[gray]{#1}}      \expandafter\def\csname LTw\endcsname{\color{white}}      \expandafter\def\csname LTb\endcsname{\color{black}}      \expandafter\def\csname LTa\endcsname{\color{black}}      \expandafter\def\csname LT0\endcsname{\color{black}}      \expandafter\def\csname LT1\endcsname{\color{black}}      \expandafter\def\csname LT2\endcsname{\color{black}}      \expandafter\def\csname LT3\endcsname{\color{black}}      \expandafter\def\csname LT4\endcsname{\color{black}}      \expandafter\def\csname LT5\endcsname{\color{black}}      \expandafter\def\csname LT6\endcsname{\color{black}}      \expandafter\def\csname LT7\endcsname{\color{black}}      \expandafter\def\csname LT8\endcsname{\color{black}}    \fi
  \fi
  \setlength{\unitlength}{0.0500bp}  \begin{picture}(7200.00,5040.00)    \gplgaddtomacro\gplbacktext{      \csname LTb\endcsname      \put(176,2574){\rotatebox{-270}{\makebox(0,0){\strut{}$F(x)=  x + 0.3x^2\sin(25x)$}}}      \put(3599,154){\makebox(0,0){\strut{}$X = [-1, 1]$}}      \put(3599,4445){\makebox(0,0){\strut{}Oriented but non-monotone operator $F(x)$}}      \put(2959,3046){\makebox(0,0)[l]{\strut{}$F(x^\star)=0$}}    }    \gplgaddtomacro\gplfronttext{    }    \gplbacktext
    \put(0,0){\includegraphics{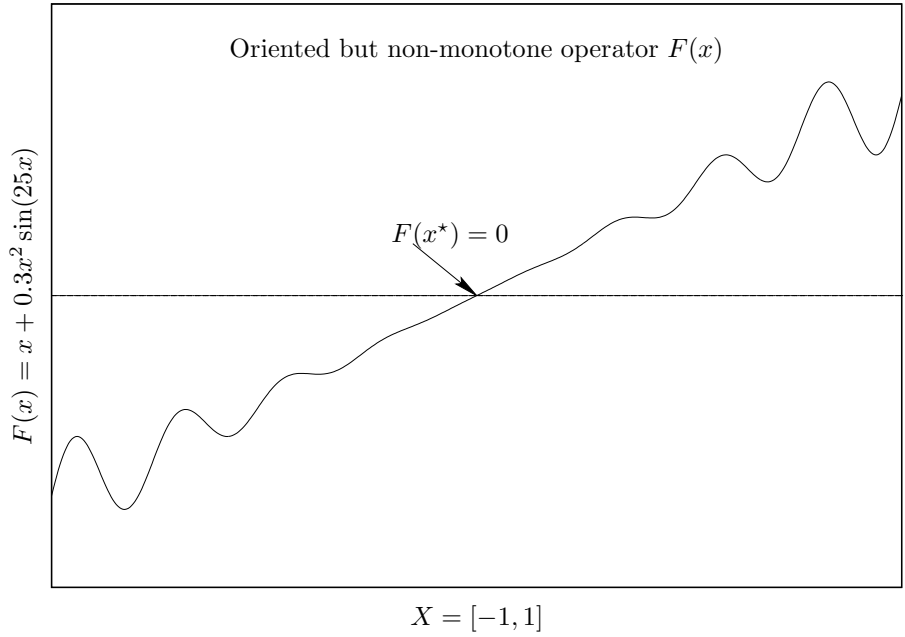}}    \gplfronttext
  \end{picture}\endgroup
 
\caption{Non-monotone operator \(F(x)\) oriented toward \(\xs = 0 \).}
\label{mono-ori}
\end{center}
\end{figure}

If \(F\) is monotone, then the pseudo-variational inequality (PVI) problem
\bleq{psvip}{ \mbox{Find } \xs \in X 
\mbox{ such that } F(x)(x - \xs) \geq 0 \mbox{ for all } x \in X.}
has a solution \(\xs\) which is  a solution of (\ref{vip}) as well.
However it is not necessary for \(F\) to be monotone to have a solution of (\ref{psvip})
which coincides with a solution of (\ref{vip})
as Fig.  \ref{mono-ori} demonstrates.

For simplicity we assume that both problems (\ref{vip}) and (\ref{psvip})
has unique and hence coinciding solutions.

To suggest a superiorized iteration method for the problem
(\ref{psvip})
we consider oriented and strongly oriented mappings according to the following definition.
\begin{definition}
A set-valued mapping \(G: E \to \cC(E)\) is called oriented toward \(\bar x\) at point \(x\) if
\bleq{Feb}{ g_x(x - \bar x) \geq 0}
for any \(g_x \in G(x)\).
\end{definition}
If \(G\) is oriented toward \(\bar x\) at all points \(x \in X\)
then we will call it oriented toward \(\bar x\) on \(X\).
Of course if \(\bar x = \xs\), a solution of PVI problem (\ref{psvip}),
then \(G\) is oriented toward \(\xs\) on \(X\) by definition and the other way around.

The notion of oriented mappings is somewhat related to attracting mappings introduced in
\cite{bb96},
which can be defined for our purposes as follows.
\begin{definition}
A mapping \(F: E \to E\) is called attracting with respect to \(\bar x\) at point \(x\) if
\bleq{attar}{ \| F(x) - \bar x \| \leq \| x - \bar x\|}
\end{definition}
It is easy to show that if \(F\) is an attracting mapping,
then \(G(x) = x - F(x)\) is an oriented mapping.
Indeed
\[
\begin{array}{c}
G(x)(x-\xb) = (x - F(x))(x - \xb) = ( x - \xb + \xb - F(x))( x - \xb) =
\|x - \xb\|^2 + (\xb - F(x))(x - \xb) \geq
\\
\|x - \xb\|^2 - \|\xb - F(x)\|\|x - \xb\| = \|x - \xb\|(\|x - \xb\| - \|F(x) - \xb\|) \geq 0.
\end{array}
\]
The reverse is not true,
\(G(x) = 10x\) is the oriented mapping toward \(\{0\}\) on \([-1, 1]\) but
neither \(G(x)\) nor \(F(x) = x - G(x) = -9x \) are attracting.

To simplify some future estimates we introduce also a following technical definition,
which can be relaxed in many ways.
\begin{definition}
\label{aldef}
A set-valued mapping \(G: E \to \cC(E)\) is called strongly oriented toward \(\bar x\) on a set \(X\) if
for any \(\epsilon > 0\) there is \(\gamma_\epsilon > 0\) such that
\bleq{alprop}{ g_x(x - \bar x) \geq \gamma_\epsilon }
for any \(g_x \in G(x)\) and all \( x \in X \setminus \{ \bar x + \epsilon B \}\).
\end{definition}

Despite the fact that the problem (\ref{vip}) depends upon the behavior of \(F\) on \(X\) only,
we need to make an additional assumption about global properties of \(F\)
to avoid certain problems with possible divergence of iteration method
due to ''run-away'' effect.
Such assumption is the long-range orientation of \(F\) which is frequently used
to ensure the desirable global behavior of iteration methods. 
\begin{definition}
A mapping \(F: E \to E\) is called long-range oriented toward a set \(X\)
if there exists
\(\rho_F \geq 0\) and \(\kappa > 0\) such that
for any \(\xb \in X \)
\bleq{coer}{ F(x)(x - \bar x) > \kappa\| x - \xb \| \mbox{ for all } x \mbox{ such that } \rho(x, X) \geq \rho_F }.
\end{definition}
We will call \(\rho_F\) the radius of longe-range orientation of \(F\) toward \(X\).
\section{Sharp penalty mapping}
\label{exap}
In this section we present the key construction
which makes possible to reduce an approximate solution of
VI problem into calculation of the limit points of iterative process,
governed by strongly oriented operators.

For this purpose we
modify
slightly
the classical definition of a polar cone of a set \(X\).
\begin{definition}
The set \( K_X(x) = \{ p: p(x-y) \geq 0 \mbox{ for all } y \in X \} \) we will call
the polar cone of \(X\) at a point \(x\).
\end{definition}
Of course \(K_X(x) = \{ 0 \}\) if \(x \in \intr X\).

For our purposes we need also a stronger definition which defines a certain sub-cone of \(K_X(x)\)
with stronger pointing toward \(X\).
\begin{definition}
Let \(\epsilon \geq 0\) and \(x \notin X + \epsilon B\).
The set
\bleq{epsicon}{
K_X^{\epsilon}(x) = \{ p: p(x-y) \geq 0 \mbox{ for all } y \in X+\epsilon B \}
}
will be called \(\epsilon\)-strong polar cone of \(X\) at \(x\).
\end{definition}
As it is easy to see that the alternative definition of \( K_X^{\epsilon}(x) \) is
\( K_X^{\epsilon}(x) = \{ p: p(x-y) \geq \epsilon \|p\| \mbox{ for all } y \in X. \} \)

To define a sharp penalty mapping for the whole space \(E\) we introduce
a composite mapping
\[
\tilde K_X^\epsilon(x) = \left\{
\begin{array}{ll}
\{ 0 \} & \mbox{ if } x \in X
\\
K_X(x) & \mbox{ if } x \in \cl \{\{ X + \epsilon B \} \setminus X \} 
\\
K_X^{\epsilon}(x) & \mbox{ if } x \in \rho_F B \setminus \{ X + \epsilon B \}
\end{array}
\right.
\]
Notice that \( \tilde K_X^\epsilon(x) \) is \usc by construction.

Now we define an \(\epsilon\)-sharp penalty mapping for \(X\) as
\[
P_X^\epsilon(x) = \{ p \in \tilde K_X^{\epsilon}(x) , \|p\| = 1 \}.
\]
Clear that \( P_X^\epsilon(x) \) is not defined for \( x \in \intr \{X\} \)
but we can defined it to be equal to zero on \(\intr \{X\} \) and take
a convex envelope of \( P_X^\epsilon(x) \) and \(\{0\}\) at the boundary of \(X\)
to preserve upper semi-continuity.

For some positive \(\lambda\) define
\( F_\lambda(x) = F(x) + \lambda P_X^\epsilon (x).\)
Of course by construction \(F_\lambda(x)\) is upper semi-continuous
for \(x \notin X\).

For the further development we establish the following result on construction
of an approximate globally oriented mapping related to the VI problem (\ref{vip}). 
\begin{enlem}
\label{lemref}
Let \(X \subset E\) is closed and bounded,
\(F\) is monotone and longe-range oriented toward \(X\) with the radius of orientability \(\rho_F\) and
strongly oriented toward solution \(\xs\) of (\ref{vip}) on \(X\)
with the constants \(\gamma_\epsilon > 0\) for \(\epsilon > 0\),
satisfying (\ref{alprop}) and
\(P_X^\epsilon(\cdot)\) is a sharp penalty.
Then for any sufficiently small \(\epsilon > 0\) there exists \(\lambda_\epsilon > 0\) and 
\(\delta_\epsilon > 0\) such that
for all \(\lambda \geq \lambda_\epsilon \) a penalized mapping
\( F_\lambda(x) = F(x) + \lambda P_X^\epsilon(x)\) satisfies the inequality
\bleq{penny}{ f_x(x - \xs) \geq \delta_\epsilon }
for all \(x \in \rho_F B \setminus \{ \xs + \epsilon B \}\) and any \(f_x \in F_\lambda(x)\).
\end{enlem}
\enProof\xspace
For monotone \(F\) we can equivalently consider a pseudo-variational inequality (\ref{psvip})
with the same solution \(\xs\).
Define the following subsets of \(E\):
\[
\begin{array}{l}
X_\epsilon^{(1)} =
X \setminus \{\xs + \epsilon B \},
\\
X_\epsilon^{(2)} =
\{
\{ X + \epsilon B \}\setminus X \} \setminus \{\xs + \epsilon B \},
\\
X_\epsilon^{(3)} =
\rho_F B \setminus
\{
\{ X + \epsilon B \} \setminus \{\xs + \epsilon B \}
\}.
\end{array}
\]
Correspondingly we consider 3 cases.
\paragraph{Case A. \( x \in X_\epsilon^{(1)}\).}
In this case \( f_\lambda(x) = F(x) \) and therefore
\bleq{caseA}{
f_\lambda(x)(x - \xs) = F(x)(x - \xs) \geq \gamma_\epsilon > 0.
}
\paragraph{Case B. \( x \in X_\epsilon^{(2)}\).}
In this case \( f_\lambda(x) = F(x) + \lambda p_X(x) \) where
\( p_X(x) \in K_X(x),\ \|p_X(x)\| = 1 \)  and therefore
\bleq{caseB}{
f_\lambda(x)(x - \xs) = F(x)(x - \xs ) + \lambda p_X(x) ( x - \xs) \geq 
\gamma_\epsilon/2 > 0.
}
as \( \lambda p_X(x) ( x - \xs) > 0 \) by construction.
\paragraph{Case C. \( x \in X_\epsilon^{(3)}\).}
In this case \( f_\lambda(x) = F(x) + \lambda p_X(x) \) where
\( p_X(x) \in K_X^\epsilon(x),\ \|p_X(x)\| = 1 \).
By continuity of \(F\)
the norm of \(F\) is bounded
on \(\rho_F B\)
by some \(M\) and as \(P^\epsilon_X(\cdot)\) is \(\epsilon\)-strong penalty mapping
\bleq{caseC}{
\begin{array}{c}
f_\lambda(x)(x - \xs) = F(x)(x - \xs ) + \lambda p_X(x) ( x - \xs) \geq 
\\
-M \| x - \xs \| + \lambda \epsilon \leq -2 \rho_F M + \lambda \epsilon \geq
\rho_F M > 0
\end{array}
}
for \(\lambda \geq \rho_F M/\epsilon \).

By combining all three bounds we obtain
\bleq{allABC}{
f_\lambda(x)(x - \xs) \geq \min \{ \gamma_\epsilon/2, \rho_F M \} = \delta_\epsilon > 0
}
for \(\lambda \geq \Lambda_\epsilon = \rho_F M/\epsilon \)
which  completes the proof.
\EndProof

The elements of a polar cone for a given set \(X\) can be obtained by different means.
The most common are either by projection onto set \(X\):
\[
x - \Pi_X(x) \in K_X(x)
\]
where \(\Pi_X(x) \in X \) is the orthogonal projection of \(x\) on \(X\),
or by subdifferential calculus when \(X\) is described by a convex inequality
\(X = \{ x: h(x) \leq 0 \}\).
If there is a point \(\xb \) such that \(h(\xb) < 0\) ( Slater condition)
then \(h(y) < 0 \) for all \( y \in \intr\{X\}\).
Therefore
\( 0 < h(x) - h(y) \leq g_h(x)(x-y) \) for any \(y \in \intr\{X\}\).
By continuity \( 0 < h(x) - h(y) \leq g_h(x)(x-y) \) for all \(y \in X\)
which means that \(g_h \in K_X(x)\).

One more way to obtain \(g_h \in K_X(x)\) relies on the ability to find some
\( x^c \in \intr \{X\} \) and use it to compute Minkowski function
\[
\mu_X(x, x^c) = \inf_{\theta \geq 0} \{ \theta: x^c + (x - x^c)\theta^{-1} \in X \}.
\]
Then 
\( \xb = x^c + (x - x^c)\mu_X(x,x^c)^{-1} \in \partial X \),
i.e. \(h(\xb) = 0\) and for any \(g_h \in \partial h(\xb) \)
holds
\(g_h \xb  \geq g_h y \) for any \(y \in X\),
in particular
\(g_h \xb \geq g_h x^c\).
Given definition of \(\xb\) obtain
\[
g_h \xb = g_h x^c + g_h (x - x^c) \mu_X(x,x^c)^{-1} =
g_h x \mu_X(x,x^c)^{-1} + ( 1 - \mu_X(x,x^c)^{-1}) g_h x^c
\]
As
\(g_h \xb \geq g_h x^c\)
we can turn equality into inequality
\[
g_h \xb \leq g_h x \mu_X(x,x^c)^{-1} + ( 1 - \mu_X(x,x^c)^{-1}) g_h \xb
\]
which leads to
\( g_h x \geq g_h \xb \geq g_h y \mbox{ for any } y \in X \),
that is \(g_h \in K_X(x)\).

As for \(\epsilon\)-expansion of \(X\) it can be approximated from above (included into)
by the relaxed inequality \( X + \epsilon B \subset \{ x:  h(x) \leq L\epsilon \} \)
 where \(L\) is a Lipschitz constant in an appropriate
neighborhood of \(X\).

\section{Convergence theory}
\makeatletter{}\label{cct}
To study the convergence we use the convergence conditions developed in \cite{nurmi-rus91}.
Within this framework the algorithm for solving a particular problem is considered as a rule
for construction of a sequence of approximate solutions
\(\{x^k\}\),
which has to converge to a certain set
\(X_\star\)
which by definition is a set of desirable solutions.
Typically this set is defined by optimality conditions (for optimization problems),
and the like.

The weak form of convergence (existence of converging sub-sequence) is guaranteed
if the sequence \(\{x^k\}\) has the following properties:
\begin{enumerate}
\item[{\bf A1}]
The sequence \(\{x^k\}\) is bounded.
\item[{\bf A2}]
There exists continuous function 
\(W(x): E \to \mathbb R\) such that if
\(\{x^k\}\) has a limit point \( x' \notin X_\star \)
then this sequence has another limit point
\(x''\)
such that \(W(x'') < W(x')\).
\end{enumerate}
If these requirements are satisfied then
the sequence \(\{x^k\}\) has a limit point
\(x^\star \in X_\star.\)

This statement is practically obvious: consider a sub-sequence
\(\{ x^{k_t}, t = 0,1, \dots \}\) such that
\[
\lim_{t \to \infty} W(x^{k_t}) = 
\lim_{n \to \infty} \inf_{m \geq n} W(x^m) = W_\star > -\infty
\]
due to continuity of \(W\) on any bounded closed set, containing \(\{ x^{k}, kt = 0,1, \dots \}\).
Then any limit point of
\(\{ x^{k_t} \}\) belongs to \(X_\star\)
otherwise 
using {\bf A2} we arrive to contradiction. 

In practice however the function \(W(\cdot)\) may be defined implicitly so the problem of selecting
the desired sub-sequences may be not so simple.
 
\section{Algorithms}
\makeatletter{}\label{salgo}
After construction of the mapping \(F_\lambda\), oriented toward solution
\(\xs\) of (\ref{psvip}) at the whole space \(E\) except \(\epsilon\)-neighborhood of \(\xs\)
we can use it in an iterative manner like
\bleq{fejit}{ x^{k+1} = x^k - \theta_k f^k,\ f^k \in F_\lambda(x^k), \ k = 0,1, \dots, }
where \(\{\theta_k\}\) is a certain prescribed sequence of step-size multipliers,
to get the sequence of \(\{x^k\}, k = 0,1, \dots \) which hopefully converges under some conditions to 
an approximate solution of (\ref{vip}).

As we use conditions {\bf A1, A2} to check for convergence we need first to establish
boundness of \(\{x^k\}, k = 0,1, \dots \).
The simplest way to guarantee this is to 
insert into the simple scheme (\ref{fejit}) a safety device
which enforces restart if a current iteration \(x^k\) goes too far.
This prevents the algorithm from  divergence due to ''run away'' effect and
keeps a sequence of iterations \(\{x^k\}\) bounded.

Than the final form of the algorithm looks like following,
where we assume that the set \(X\), VI operator \(F\) and sharp penalty mapping \(P_X\) satisfy
conditions of the lemma \ref{lemref}:
\begin{center}
\begin{algorithm}[H]
\SetAlgoLined
\KwData{The variational inequality operator \(F\), sharp penalty mapping \(P_X\),
positive constant \(\epsilon\), 
penalty constant \(\lambda > \Lambda_\epsilon,\) which existence is claimed by the Lemma \ref{lemref},
longe-range orientation radius \(\rho_F\),
a sequence of step-size multipliers \(\{0 < \theta_k,\, k = 0,1,2, \dots \}\).
and an initial point \(x^0 \in \rho_F B\).
}

\KwResult{The sequence of approximate solutions \(\{x^k\}\) which contains
a converging sub-sequence \(\{ x^{t_k}\}\) the limit point of which belongs
to an \(\epsilon\)-solution of variational inequality (\ref{vip}).}

{\bf Initialization}\;
Define penalized mapping
\[
F_\lambda(x) = F(x) + \lambda P_X(x),
\]
and set the iteration counter \(k\) to \(0\)\;

\While { The limit is not reached }{
Generate a next approximate solution \( x_{k+1} \):
\bleq{algo}{
x^{k+1} = \left\lbrace
\begin{array}{ll}
x^k - \theta_k f^k,\ f^k \in F_\lambda(x^k),
& \mbox{ if } \|x^k\| \leq 2\rho_F
\\
x^0	& \mbox{ otherwise.}
\end{array}
\right.
}
Increment iteration counter \(k \longrightarrow k+1 \)\;
}
{\bf Complete:} accept
\(\{x^k\}, k = 0,1, \dots \)
as a solution of (\ref{vip})
\footnote{
This is to mean that we still face two problems: select a sub-sequence
which converges to an \(\epsilon\)-solution and provide a stopping criteria for that.}
\caption{The generic structure of the fixed point iteration algorithm with exact penalty.}
\label{algosub}
\end{algorithm}
\end{center}
Of course this is a conceptual version of the algorithm as it has no termination criteria.

By rewriting (\ref{algo}) as
\bleq{algo2}{
x^{k+1} = \left \{
\begin{array}{ll}
x^k - \lambda\theta_k \lambda^{-1}f^k,\ f^k \in F_\lambda(x^k) & \mbox{ if } \|x^k\| \leq 2\rho_F
\\
x^0	& \mbox{ otherwise.}
\end{array}
\right.
}
and redefining \(\lambda\theta_k \to \theta_k \) and
\( \lambda^{-1}F_\lambda(x) = P_X(x) + \lambda^{-1}F(x) \to F_\lambda(x) \)
we may see the effect of superiorization of the feasibility iterative algorithm
\( x^{k+1} = x^k + \theta_k p^k,\ p^k \in P_X(x^k), k = 0,1,\dots \) toward the algorithm (\ref{algo2})
for solution of the variational inequality (\ref{vip}) with the help
of small perturbation \(\lambda^{-1}F(x^k)\).
The effect of these perturbations when \(\lambda^{-1} \to 0\) was discovered in general case
in \cite{nurmiDAN08,nurmiCMMP08} here we see that it can be achieved with a finite \(\lambda\).

We prove convergence of the algorithm \ref{algosub} under common assumptions on step sizes
\( \theta_k \to +0\) when \(k \to \infty\) and \(\sum_{k=1}^K \theta_k \to \infty\) when \(K\to\infty\).
This is not the most efficient way to control the algorithm, but at the moment we are interested mostly in
the very fact of convergence.
\begin{theorem}
\label{thmref}
Let \(\epsilon > 0, F,  P_X^\epsilon\) satisfy the assumptions of the Lemma \ref{lemref},
and \(\Lambda_\epsilon\) is such that \(F_\lambda(x) = F(x) + \lambda P_X^\epsilon(x)\) is
\(\delta_\epsilon\)-oriented with respect to \(\xs\) on
\(\rho_F B \setminus \{\xs + \epsilon B \} \) with \(\delta_\epsilon > 0\)
for any \(\lambda \geq \Lambda_\epsilon\),
and \( \theta_k \to +0\) when \(k \to \infty\) and \(\sum_{k=1}^K \theta_k \to \infty\) when \(K\to\infty\).
Then there is a limit point of the sequence \(\{x^k\}\) generated by the algorithm \ref{algosub}
which belongs to the set of \(\epsilon\)-solutions \(\xs + \epsilon B \) of the problem (\ref{vip}). 
\end{theorem}
\enProof\xspace
Notice first that \(\Lambda_\epsilon\) and \(\delta_epsilon > 0\)
exist due to lemma \ref{lemref}.

We prove the theorem by demonstrating that
the Algorithm 1 satisfies convergence conditions {\bf A1,A2}
of Section \ref{cct}.
The most basic property asked for in these conditions is the boundness
of the algorithm trajectory \(\{x^k, k = 1,2, \dots\}\) and we show it first.

Let \(\bar\rho_F\) be the radius of long-range orientability of the operator \(F\) and
\(\rho_F\) in Algorithm 1 is large enough that
\(X + \bar\rho_F B \subset X + 3 \bar\rho_F/2 B \subset \rho_F B\).

If the sequence \(\{\|x^k\|\}\) leaves \(\rho_F B\) infinitely many times
then it should leave the set \(X + 3 \bar\rho_F/2 B \) infinitely many times as well. 
Define ( a finite or not ) set \(T\) of indexes \(T = \{t_k, k = 1,2, \dots\}\) such that
\[
\| x^{\tau} - \xs \| < \frac32 \bar\rho_F \mbox{ and } \| x^{\tau + 1} - \xs \| \geq \frac32 \bar\rho_F.
\]
and denote for brevity \(f^\tau = f_\lambda(x^\tau) = F(x^\tau) + \lambda p_X^\epsilon(x^\tau),\)
where \(p_X^\epsilon(x^\tau) \in P_X^\epsilon(x^\tau)\).
If \(\tau \in T\) then
\[
\| x^{\tau + 1} - \xs \|^2 = \| x^{\tau}  - \theta_k f^\tau -\xs \|^2 =
\| x^{\tau}  - \xs \|^2 - 2 \theta_\tau f^\tau (x^{\tau} -\xs ) + \theta_\tau^2 \| f^{\tau}\|^2 
\]
Due to Lemma \ref{lemref} \(f_\lambda(\cdot)\) is oriented toward \(\xs\)
therefore \( f^\tau (x^{\tau} -\xs ) \geq \delta_\epsilon > 0 \) and due to
upper semi-continuity of \(f_\lambda(\cdot)\)
there is \(C\) such that \(\| f^{\tau}\|^2 \leq C.\)
Hence
\[
\| x^{\tau + 1} - \xs \|^2 \leq \| x^{\tau}  - \xs \|^2 - 2 \theta_\tau \delta_\epsilon + \theta_\tau^2 C
\]
and as \( 2 \delta_\epsilon - \theta_\tau C > \delta_\epsilon \) 
for large enough \(\tau\)
\[
\| x^{\tau + 1} - \xs \|^2 \leq \| x^{\tau}  - \xs \|^2 - \theta_\tau \delta_\epsilon < \| x^{\tau}  - \xs \|^2
\]
which is a contradiction with a choice of \(\tau \in T\).
Therefore \(T\) is a finite set and
the sequence \(\{x^k\}\) leaves the set \(\rho_F B\) a finite number of times only.
It proves the condition {\bf A1} of convergence conditions of the section \ref{cct}.

The next step is to prove that the other condition {\bf A2} is fulfilled as well.
For this we define \(W(x) = \| x - \xs \|^2\) and demonstrate that it is indeed so.

If \(x^{t_k} \to x' \notin X_\star\) then \(f'(x' - \xs) > 0 \) for any \(f' \in F_\lambda(x')\)
and by upper semi-continuity of \(F_\lambda\) there exists an \(\varepsilon > 0\)
such that \(F_\lambda(x)(x - \xs) \geq \delta\) for all \(x \in x'+ 4\varepsilon B\) and some \(\delta > 0\).

For for \(k\) large enough \( x^{t_k} \in x' + \varepsilon B\) and let us assume that for all \( t > t_k \)
the sequence \(\{x^t, t > t_k \} \subset x^{t_k} + \varepsilon B \subset x' + 2\varepsilon B\).
Then again
\bleq{exit}{
\begin{array}{c}
\| x^{t + 1} - \xs \|^2 =
\| x^{t}  - \theta_t F_\lambda(x^t) -\xs \|^2 =
\| x^t  - \xs \|^2 - 2 \theta_t F_\lambda(x^t)(x^t -\xs ) + \theta_t^2 \| F_\lambda(x^t)\|^2 \leq
\\
\| x^t  - \xs \|^2 - 2 \theta_t F_\lambda(x^t)(x^t -\xs ) + \theta_t^2 C^2 \leq
\| x^t  - \xs \|^2 - 2 \theta_t \delta + \theta_t^2 C^2 <
\| x^t  - \xs \|^2 - \theta_t \delta,
\end{array}
}
for all \( t > t_k \) and \(k\) large enough that \( \sup_{t > t_k} \theta_t < \delta/C^2 \). 
Summing up last inequalities from \(t =t_k\) to \(t = T-1\) obtain
\bleq{lyap}{
\| x^T - \xs \|^2 \leq \| x^{t_k} - \xs \|^2 - \delta \sum_{t = t_k}^{T-1} \theta_t \to -\infty
}
when \(T \to \infty\) which is of course impossible.

Hence for each \(t_k\) there exists \(s_k > t_k\) such that \(\|x^{t_k} - x^{s_k} \| > \varepsilon > 0 \)
Assume that \(s_k\) in fact a minimal such index, i.e.
\( \| x^{t_k} - x^t \| \leq \varepsilon  \) for all \(t\) such that \( t_k < t < s_k \)
that is \(x^t \in x^{t_k} + \varepsilon B \subset x' + 2\varepsilon B\).
Without any loss of generality we can assume \(x^{s_k} \to x''\) where by construction
\(\| x' - x''\| \geq \varepsilon > 0\) and therefore \(x' \neq x''\).

As all conditions which led to (\ref{lyap}) hold for \(T = s_k\) then by letting \(T= s_k\) obtain
\[
W(x^{s_k}) \leq W(x^{t_k}) - \delta \sum_{t = t_k}^{s_k-1} \theta_t.
\]
On the other hand
\[
\varepsilon < \|x^{t_k} - x^{s_k}\| \leq \sum_{t = t_k}^{s_k-1} \| x^{t+1} - x^t \|
\leq \sum_{t = t_k}^{s_k-1} \theta_t \|F_\lambda(x^t)\| \leq K \sum_{t = t_k}^{s_k-1} \theta_t
\]
where \(K\) is the upper estimate of the norm of \(F_\lambda(x)\) on \(2\rho_F B\).

Therefore \(\sum_{t = t_k}^{s_k-1} \theta_t > \varepsilon/K > 0 \) and finally
\[
W(x^{s_k}) \leq W(x^{t_k}) - \delta\varepsilon/K.
\]
Passing to the limit when \(k\to\infty\) obtain \(W(x'') \leq W(x') - \delta \varepsilon/K < W(x') \)
which proves {\bf A2} and therefore completes the proof.
\EndProof
 
\section*{Conclusions}
In this paper we made use of a sharp penalty mapping to construct iteration algorithm
weakly converging to an approximate solution of a monotone variational inequality.
It can be considered as a variant of a superiorization technique which combines feasibility
and optimization steps into joint process but applied to the different type of problems.
These problems create however an additional theoretical difficulty when making them amendable to simple iteration algorithm
and it was possible in this paper to prove only weak approximate convergence result (convergence of sub-sequence to
\(\epsilon\)-solution).

As for practical value of these result
it is generally believed that the conditions for the step-size multipliers used in this theorem result
in rather slow convergence of the order \(O(k^{-1})\).
However the convergence rate can be improved by different means following the example 
of non-differentiable optimization.
The promising direction is f.i. the least-norm adaptive regulation,
suggested probably first by A.Fiacco and McCormick
\cite{sumt68}
as early as 1968
and studied in more details in
\cite{oms}
for convex optimization problems.
With some modification in can be easily used for VI problems as well.
Experiments show that under favorable conditions it produces step multipliers decreasing
as geometrical progression which gives a linear convergence for the algorithm.
This may explain the success of 
\cite{bdhk2007}
where geometrical progression for step multipliers was independently suggested and tested in practice.

\end{document}